\documentstyle[12pt]{article}
\begin{document}
\newpage
\pagestyle{empty}
\setcounter{page}{0}
%%%%%%%%%%%%%
\renewcommand{\theequation}{\thesection.\arabic{equation}}
\newcommand{\sect}[1]{\setcounter{equation}{0}\section{#1}}
%%%%%%%%%%%%%
%%% ** start of amsfont definitions **
\newfont{\twelvemsb}{msbm10 scaled\magstep1}
\newfont{\eightmsb}{msbm8}
\newfont{\sixmsb}{msbm6}
\newfam\msbfam
\textfont\msbfam=\twelvemsb
\scriptfont\msbfam=\eightmsb
\scriptscriptfont\msbfam=\sixmsb
\catcode`\@=11
\def\Bbb{\ifmmode\let\next\Bbb@\else
  \def\next{\errmessage{Use \string\Bbb\space only in math mode}}\fi\next}
\def\Bbb@#1{{\Bbb@@{#1}}}
\def\Bbb@@#1{\fam\msbfam#1}
\newfont{\twelvegoth}{eufm10 scaled\magstep1}
\newfont{\tengoth}{eufm10}
\newfont{\eightgoth}{eufm8}
\newfont{\sixgoth}{eufm6}
\newfam\gothfam
\textfont\gothfam=\twelvegoth
\scriptfont\gothfam=\eightgoth
\scriptscriptfont\gothfam=\sixgoth
\def\frak{\frak@}
\def\frak@#1{{\fam\gothfam{{#1}}}}
\def\frak@@#1{\fam\gothfam#1}
\catcode`@=12
%%% ** end of amsfont definitions **
%\def\Bbb{\bf}
\def\arcsinh{\mathop{\rm arcsinh}\nolimits}
\def\CC{{\Bbb C}}
\def\NN{{\Bbb N}}
\def\QQ{{\Bbb Q}}
\def\RR{{\Bbb R}}
\def\ZZ{{\Bbb Z}}
\def\cA{{\cal A}}          \def\cB{{\cal B}}          \def\cC{{\cal C}}
\def\cD{{\cal D}}          \def\cE{{\cal E}}          \def\cF{{\cal F}}
\def\cG{{\cal G}}          \def\cH{{\cal H}}          \def\cI{{\cal I}}
\def\cJ{{\cal J}}          \def\cK{{\cal K}}          \def\cL{{\cal L}} 
\def\cM{{\cal M}}          \def\cN{{\cal N}}          \def\cO{{\cal O}}
\def\cP{{\cal P}}          \def\cQ{{\cal Q}}          \def\cR{{\cal R}} 
\def\cS{{\cal S}}          \def\cT{{\cal T}}          \def\cU{{\cal U}}
\def\cV{{\cal V}}          \def\cW{{\cal W}}          \def\cX{{\cal X}}
\def\cY{{\cal Y}}          \def\cZ{{\cal Z}}
\def\id{\mbox{id}}
\def\ggo{{\frak g}_{\bar 0}}
\def\uqggo{\cU_q({\frak g}_{\bar 0})}
\def\uqggp{\cU_q({\frak g}_+)}
\def\typeA{{\em type $\cA$}}
\def\typeB{{\em type $\cB$}}

\newtheorem{Def}{Definition}
\newtheorem{lemme}{Lemme}
\newtheorem{prop}{Proposition}
\newtheorem{theo}{Th\'eor\`eme}
\vfill

$$
\;
$$

\vskip 2.5cm

\begin{center}

  {\LARGE {\bf {\sf On a Deformation of $sl(2)$ with Paragrassmannian 
Variables}}} \\[1cm]

\vfill

{\large B. Abdesselam$^{\dagger,2}$,
  J. Beckers$^{\ddagger,3}$, A. Chakrabarti$^{\dagger,4}$
  and N. Debergh$^{\ddagger,3,5}$}

\vfill

{\em $^{\dagger}$
  Centre de Physique Th{\'e}orique
  \footnote{Laboratoire Propre du CNRS UPR A.0014  

  \indent
  $\;$$^2$abdess@orphee.polytechnique.fr,

  \indent
  $\;$$^3$beckers@vm1.ulg.ac.be,

  \indent
  $\;$$^4$chakra@orphee.polytechnique.fr,

  \indent
  $\;$$^5$Chercheur Institut Interuniversitaire des Sciences 
  Nucl{\'e}aires (Brussels).}, Ecole Polytechnique, 
  91128 Palaiseau Cedex,
  France.
  \indent
  \vskip 0.2cm
  $^{\ddagger}$ Physique Th{\'e}orique et Math{\'e}matique, Institut de 
  Physique, B5, Universit{\'e} de Li{\`e}ge, B-4000-Li{\` e}ge 1, Belgium.}

\end{center}

\vfill

\begin{abstract}
We propose a new structure ${\cal U}^{r}_{\displaystyle{q}}(sl(2))$. 
This is realized by multiplying $\delta$ ($q=e^{\delta}$, $\delta\in \CC$) 
by $\theta$, where $\theta$ is a real nilpotent -paragrassmannian- variable 
of order $r$ ($\theta^{r+1}=0$) that 
we call the order of deformation, the limit $r\rightarrow \infty$ giving 
back the standard ${\cal U}_{\displaystyle{q}}(sl(2))$. In particular
we show that, for $r=1$, there exists a new ${\cal R}$-matrix associated 
with $sl(2)$. We also proof that the
restriction of the values of the parameters of deformation give 
nonlinear algebras as particular cases. 
\end{abstract}

\vfill
\vfill
\vfill
\vfill

\newpage
\pagestyle{plain}

\sect{Introduction}

During the last few years, $q$-deformations $[1]$ ($q=e^{\delta}$, 
$\delta \in \CC$) of the universal enveloping algebra of Lie algebras have 
attracted a wide attention. They are indeed remarkable mathematical 
structures known as Hopf algebras and they have been proved to be connected 
to Conformal Field Theory, in particular, as they have been figuring in 
$2d$-solvable model $S$-matrices and solutions to their Yang-Baxter 
factorization equations (See Ref. $[2]$ and references therein).   

\smallskip

The pioneering papers $[3]$ devoted to the specific 
${\cal U}_{\displaystyle{q}}(sl(2))$ case have been extended by various 
authors. Let us just mention here the Ro${\tilde c}$ek proposal $[4]$ (based 
on generalized nonlinear deformations) providing a new algebraic 
description of the Morse and modified P$\ddot o$schl-Teller Hamiltonians 
$[5]$. Despite of its physical interest, the Ro${\tilde c}$ek deformation 
has been rarely exploited, compared to the Drinfeld-Jimbo one, because of 
its mathematical defect: its Hopf characteristics (coproduct, counit, 
antipode) have not yet been pointed out. 

\smallskip

In this letter, we answer the following question: Is it possible to obtain 
the nonlinear algebras as particular restrictions of the quantum deformation ?

\smallskip

Our purpose is then twofold. First, we introduce the nilotent algebra
${\cal U}^{r}_{\displaystyle{q}}(sl(2))$ by multiplying $\delta$ by $\theta$,
 where $\theta$ is a real nilotent -paragrassmannian- 
variable $[6]$ of order $r$ ($\theta^{r+1}=0$). Second, we discuss the 
connection of this new structure to some particular nonlinear deformations 
of $sl(2)$ whose Hopf characteristics are introduced.  

\smallskip 

In section 2, we briefly review the Drinfeld-Jimbo deformation of $sl(2)$. 
Then, in section 3, we introduce the quantization with one paragrassmannian 
variable and its Hopf structure. The 
quantization with two paragrassmannian variables is given in section 4. In 
section 5, we give the connection of these structures to particular 
nonlinear deformations of $sl(2)$. Finally, we conclude in section 6 with some comments.  

\sect{The ${\cal U}_{\displaystyle{q}}(sl(2))$ algebra.}

The standard Drinfeld-Jimbo deformation $[1]$ of the Lie algebra $sl(2)$ 
generated by $H$, $J_{+}$, $J_{-}$ is characterized by the relations 
\begin{eqnarray}
&&[\;J_{+},\;J_{-}]={q^{H}-q^{-H} \over q-q^{-1}}=
{\sinh (\delta H) \over \sinh (\delta)}, \nonumber\\
&&[\;H,\;J_{\pm}]=\pm 2 J_{\pm}.
\end{eqnarray}
It is completed by the additional operations, coproduct  
$\triangle :{{\cal U}_{q}}(sl(2)) \rightarrow {{\cal U}_{q}}(sl(2))
\otimes {{\cal U}_{q}}(sl(2))$, counit $\varepsilon : {{\cal U}_{q}}(sl(2))
\rightarrow\CC$ and the antipode $S:{{\cal U}_{q}}(sl(2))\rightarrow 
{{\cal U}_{q}}(sl(2))$ such that
\begin{eqnarray}
&&\bigtriangleup (H) = H \otimes \hbox{{\bf 1}} + \hbox{{\bf 1}} 
\otimes H,\nonumber \\
&&\bigtriangleup (J_{\pm}) = J_{+}\otimes e^{\delta H/2} + 
e^{-\delta H/2}\otimes J_{\pm}, \nonumber\\
&& \varepsilon (\hbox{{\bf 1}}) =\hbox{{\bf 1}}, \;\; 
\varepsilon (J_{\pm})=\varepsilon (H)=0,
\nonumber \\
&& S(\hbox{{\bf 1}})=\hbox{{\bf 1}},\;\;S(H)=-H,\;\;
S(J_{\pm})=-e^{{\pm} \delta}\;J_{\pm} 
\end{eqnarray}
where $\bigtriangleup$ and $\varepsilon$ are homomorphisms while $S$ is an 
algebra antihomomorphism
\begin{eqnarray}
&&\bigtriangleup (a\;b)= \bigtriangleup (a)\bigtriangleup (b), \nonumber\\
&&\varepsilon (a\;b)= \varepsilon (a)\varepsilon (b), \nonumber\\
&&S(a\;b)=S(b)S(a). 
\end{eqnarray}

Moreover, if $m :{{\cal U}_{q}}(sl(2))\otimes {{\cal U}_{q}}(sl(2)) \rightarrow
{{\cal U}_{q}}(sl(2))$ stands for the multiplication mapping of ${{\cal U}_{q}}(sl(2))$ i.e. $m(a \otimes b)=a.b$, we have 
\begin{eqnarray}
&&(id \otimes \bigtriangleup )\bigtriangleup =
(\bigtriangleup \otimes id)\bigtriangleup , \nonumber\\
&&m(id \otimes S)\bigtriangleup  = m(S \otimes id)\bigtriangleup =
i\circ \varepsilon , \nonumber \\
&&(\varepsilon \otimes id)\bigtriangleup = (id \otimes \varepsilon)
\bigtriangleup  =id. 
\end{eqnarray}
These are just all the axioms of a Hopf algebra, and so 
${{\cal U}_{q}}(sl(2))$ endowed with $\varepsilon$, $\bigtriangleup$ and 
$S$ just forms a Hopf algebra.

\smallskip

Let us define the formal series
\begin{eqnarray}
&&J_{\pm}=\sum_{k=0}^{\infty} \delta^{k}\;J_{\pm}^{(k)} 
\end{eqnarray}
and
\begin{eqnarray}
&&{\sinh (H \delta) \over \sinh (\delta)}=\sum_{k=0}^{\infty} \psi_{k}(H) 
\;\delta^{2k}, 
\end{eqnarray}
the second formula being just the result of a Taylor expansion. The generators 
$J_{\pm}^{(k)}$ and $H$ satisfy the following commutation relations
\begin{eqnarray}
&&[H,\;J_{\pm}^{(k)}]=\pm 2\;J_{\pm}^{(k)},\nonumber \\
&&\sum_{m=0}^{2\;k} [J_{+}^{(m)},\;J_{-}^{(2k-m)}]=\psi_{k}(H), \qquad \sum_{m=0}^{2\;k+1} [J_{+}^{(m)},\;J_{-}^{(2k+1-m)}]=0, \nonumber\\
&&\sum_{m=0}^{k}[J_{\pm}^{(m)},\;J_{\pm}^{(k-m)}]=0. 
\end{eqnarray}
Its Hopf structure is given by 
\begin{eqnarray}
&&\bigtriangleup(H)=\hbox{{\bf 1}}\otimes H+H\otimes 
\hbox{{\bf 1}},\nonumber \\
&&\bigtriangleup(J_{\pm}^{(k)})=\sum_{m=0}^{k}{1 \over 2^{m}\;m!}\bigl((-1)^{m}
H^{m}\otimes J_{\pm}^{(k-m)}+J_{\pm}^{(k-m)}\otimes H^{m}\bigr), \nonumber\\
&&\varepsilon (H)=\varepsilon (J_{\pm}^{(k)})=0,\;\varepsilon(\hbox{{\bf 1}})=
\hbox{{\bf 1}}, \nonumber\\
&&S(J_{\pm}^{(k)})=-\sum_{m=0}^{k}{(\pm)^{m} \over m!} 
J_{\pm}^{(k-m)},\;S(H)=-H,\;S(\hbox{{\bf 1}})=\hbox{{\bf 1}}, 
\end{eqnarray}
as it can be verified.

\sect{The ${\cal U}^{r}_{\displaystyle{q}}(sl(2))$ 
algebra} 

Let us introduce the real nilpotent -paragrassmannian- variable $\theta$ 
of order $r$, i.e.
\begin{eqnarray}
&\theta^{r+1}=0, 
\end{eqnarray}
being realized, in a simple way, by
\begin{eqnarray}
\theta=\pmatrix{
0&0&\cdots&0 \cr
1&0&\cdots&0& \cr
0&\ddots&\ddots&\vdots \cr
0&0&1&0 \cr}. 
\end{eqnarray}
Besides this choice, we want to notice that there are other representations 
such as that given by
\begin{eqnarray}
\theta=\sum_{\alpha=1}^{r}\theta^{(\alpha)},
\end{eqnarray}
where 
\begin{eqnarray}
(\theta^{(\alpha)})^2=0, \qquad\qquad\qquad [\theta^{(\alpha)},
\theta^{(\beta)}]=0,\;\;\;\alpha\not=\beta.   
\end{eqnarray}
Then with Eq. $(3.2)$,
we propose to generalize the operators $(2.5)$ through
\begin{eqnarray}
J_{\pm}^{\theta}&=&\sum_{m=0}^{r} \delta^{m}\;\theta^{m}\;J_{\pm}^{(m)} \\
&=&\pmatrix{
J_{\pm}^{(0)}&0&\cdots&0 \cr
\delta J_{\pm}^{(1)}&J_{\pm}^{(0)}&\cdots&0& \cr
\vdots & & & \cr
\delta^{r-1}J_{\pm}^{(r-1)}&\ddots&\ddots&\vdots \cr
\delta^{r}J_{\pm}^{(r)}&\delta^{r-1}J_{\pm}^{(r-1)}&\delta
J_{\pm}^{(1)}&J_{\pm}^{(0)} \cr}. 
\end{eqnarray}
Using the commutations relations $(2.7)$, we thus have
\begin{eqnarray}
&&[H,\;J_{\pm}^{\theta}]=\pm 2 J_{\pm}^{\theta},\qquad\qquad\qquad
\qquad\qquad\qquad
\end{eqnarray}
and
\begin{eqnarray}
&&[J_{+}^{\theta},\;J_{-}^{\theta}]=\sum_{k=0}^{r}
\delta^{k}\;\theta^{k}\;\Bigl(
\sum_{m=0}^{k}[J_{+}^{(m)},\;J_{-}^{(k-m)}]\bigr) \nonumber \\
&& \phantom{[J_{+}^{\theta},\;J_{-}^{\theta}] } =\psi_{0}(H)+
\theta^{2}\;\delta^{2}\;\psi_{1}(H)+\cdots+\theta^{2\;[r/2]}
\delta^{2\;[r/2]}\;\psi_{[r/2]}(H) \nonumber\\
&&\phantom{[J_{+}^{\theta},\;J_{-}^{\theta}]}=\sum_{k=0}^{[r/2]}
\psi_{k}(H)\;\theta^{2k}\;\delta^{2k}, 
\end{eqnarray}
where $[\lambda ]$ stands for the integer part of $\lambda$. Defining the 
exponential map by 
\begin{eqnarray}
e\;(x;\theta)=\sum_{k=0}^{r}{x^{k}\theta^{k} \over k!}, 
\end{eqnarray}
we can finally write
\begin{eqnarray}
&&[J_{+}^{\theta},\;J_{-}^{\theta}]={e\;(H\delta;\;\theta)-
 e\;(-H\delta;\;\theta) \over e\;(\delta;\;\theta)-
e\;(-\delta;\;\theta)}, \nonumber\\
&&[H,\;J_{\pm}^{\theta}]=\pm 2 J_{\pm}^{\theta}. 
\end{eqnarray}

The algebra $\lbrace J_{\pm}^{\theta},\;H\rbrace$ described by the 
commutation relations $(3.10)$ is just the deformation of $sl(2)$ with 
one paragrassmannian variable and is denoted by 
${\cal U}^{r}_{\displaystyle{q}}(sl(2))$. This algebra is
isomorphic to ${\cal U}_{\displaystyle{q}}(sl(2))/
(\delta^{r+1}{\cal U}_{\displaystyle{q}}(sl(2)))$, i.e.
\begin{eqnarray}
{\cal U}^{r}_{\displaystyle{q}}(sl(2))\cong 
{\cal U}_{\displaystyle{q}}(sl(2))/(\delta^{r+1}
{\cal U}_{\displaystyle{q}}(sl(2))). \nonumber
\end{eqnarray}

In order to define a Hopf structure for 
${\cal U}^{r}_{\displaystyle{q}}(sl(2))$, we need the following  
definition

\begin{Def} Let  
\begin{eqnarray}
a=a_{0}+a_{1} \theta+\cdots+a_{r}\theta^{r},\;\;\;\;\;\;\;b=b_{0}+
b_{1} \theta+\cdots+b_{r}\theta^{r}, 
\end{eqnarray} 
the tensor product between $a$ and $b$ is defined by
\begin{eqnarray}
a\;\bar \otimes\; b=\sum_{m=1}^{r}\sum_{n=1}^{r}a^{(m)}\otimes b^{(n)} 
\theta^{m+n}, 
\end{eqnarray}
and
\begin{eqnarray}
(a\;\bar \otimes\;b)(c\;\bar \otimes\;d)=(ac\;\bar \otimes\;bd). 
\end{eqnarray}
This operation is called the paragrassmannian tensor product. 
\end{Def}

When the paragrassmannian order $r\rightarrow \infty$, this operation is 
equivalent to the standard one. This paragrassmannian tensor 
product is compatible with 
\begin{eqnarray}
{\cal U}^{r}_{\displaystyle{q}}(sl(2))\bar 
\otimes {\cal U}^{r}_{\displaystyle{q}}(sl(2))\equiv 
{\cal U}^{r}_{\displaystyle{q}}(so(4))
\end{eqnarray}
and with the inclusion
\begin{eqnarray} {\cal U}^{r}_{\displaystyle{q}}(sl(2))\subset 
{\cal U}^{r}_{\displaystyle{q}}(sl(3))\subset
\cdots {\cal U}^{r}_{\displaystyle{q}}(sl(N-1))\subset
{\cal U}^{r}_{\displaystyle{q}}(sl(N)). 
\end{eqnarray}

We are now able to claim that

\begin{prop}
The Hopf structure associated to the ${\cal U}_{q}^{r}(sl(2))$ 
is given by
\begin{eqnarray}
&&\bigtriangleup(H)=H\otimes\hbox{{\bf 1}}+\hbox{{\bf 1}}\otimes H, 
\nonumber\\
&&\bigtriangleup(J_{\pm}^{\theta}) = J_{\pm}^{\theta}\;\bar \otimes \; 
e\;({H\delta\over 2};\;\theta)+ e\;(-{H\delta\over 2};\;\theta) 
\;\bar \otimes \; J_{\pm}^{\theta},\nonumber \\
&& \varepsilon (J_{\pm}^{\theta})=\varepsilon (H)=0,
\qquad \varepsilon(\hbox{{\bf 1}})=\hbox{{\bf 1}}, \nonumber\\
&& S(H)=-H,\qquad S(J_{\pm}^{\theta})=-e\;(\pm \delta;\;\theta)\;
J_{\pm}^{\theta}, \qquad S(\hbox{{\bf 1}})=\hbox{\bf 1},\nonumber \\
&&\bigtriangleup(e\;({H\delta\over 2};\;\theta))=
e\;({H\delta\over 2};\;\theta)\;\bar \otimes 
\; e\;({H\delta\over 2};\;\theta).
\end{eqnarray}
The following axioms are then satisfied 
\begin{eqnarray}
&&(id \;\bar \otimes \;\bigtriangleup )\bigtriangleup 
=(\bigtriangleup\; \bar \otimes \;id)\bigtriangleup ,\nonumber \\
&&m(id \;\bar \otimes\; S)\bigtriangleup 
= m(S \;\bar \otimes\; id)\bigtriangleup =
i\circ \varepsilon,\nonumber \\
&&(\varepsilon \;\bar \otimes\; id)\bigtriangleup = 
(id \;\bar \otimes \;\varepsilon)\bigtriangleup  = id, 
\end{eqnarray}
with the coproduct  
$\triangle :{{\cal U}^{r}_{q}}(sl(2)) \rightarrow 
{{\cal U}^{r}_{q}}(sl(2)) \bar 
\otimes {{\cal U}^{r}_{q}}(sl(2))$, counit $\varepsilon : 
{{\cal U}^{r}_{q}}(sl(2))$ $\rightarrow
\CC[\theta]$, the antipode $S:{{\cal U}^{r}_{q}}(sl(2))\rightarrow 
{{\cal U}^{r}_{q}}(sl(2))$ and
$m:\;{{\cal U}^{r}_{q}}(sl(2)) \bar 
\otimes {{\cal U}^{r}_{q}}(sl(2))$ $\rightarrow {{\cal U}^{r}_{q}}(sl(2))$,
where the operations $\bigtriangleup$, $S$ and $\varepsilon$ only act on $H$ 
and $J_{\pm}^{(m)}$.
\end{prop}

\smallskip

Let us now turn to some specific examples.

\smallskip

{\bf Example. 1}. The $r=0$ case is characterized by
\begin{eqnarray}
\theta =0,\;\;\;\;\;\;\;\;\;\;\;J_{\pm}^{\theta}=J_{\pm}^{(0)} \nonumber
\end{eqnarray}
and
\begin{eqnarray}
&&[H, \;J_{\pm}^{\theta}]=\pm 2 J_{\pm}^{\theta}, \nonumber\\
&&[J_{+}^{\theta},\;J_{-}^{\theta}]=H. 
\label{eq:theta=0}
\end{eqnarray}
Thus, the ${\cal U}_{q}^{0}(sl(2))$ algebra is 
nothing but $sl(2)$, endowed as usual with
\begin{eqnarray}
&&\bigtriangleup(H)=\hbox{{\bf 1}}\otimes H + H\otimes \hbox{{\bf 1}}
\nonumber\ \\
&&\bigtriangleup(J_{\pm}^{\theta})=J_{\pm}^{\theta}\otimes 
\hbox{{\bf 1}} + \hbox{{\bf 1}} \otimes 
J_{\pm}^{\theta},\;\hbox{etc.} \nonumber
\end{eqnarray} 

\smallskip

{\bf Example. 2}. The $r=1$ case is characterized by
\begin{eqnarray}
\theta=\pmatrix{
0&0 \cr
1&0 \cr}\;\;\;\;\;\;\;\hbox{and}\;\;\;\;\;\;\;J_{\pm}^{\theta}=\pmatrix{
J_{\pm}^{(0)}&0 \cr
\delta J_{\pm}^{(1)}& J_{\pm}^{(0)} \cr} \nonumber
\end{eqnarray}
and the $sl(2)$ algebra (\ref{eq:theta=0}) but now supplemented 
by a non cocommutative coproduct
\begin{eqnarray}
&&\bigtriangleup (\;H\;) = H \otimes \hbox{{\bf 1}} + 
\hbox{{\bf 1}} \otimes H,\nonumber \\
&&\bigtriangleup (\;J_{\pm}^{\theta}\;) = 
J_{\pm}^{\theta}\;\bar \otimes \; 
(\hbox{{\bf 1}}+{ \theta \delta \over 2}H) + 
(\hbox{{\bf 1}}-{\theta\delta\over 2}H)\;\bar 
\otimes \; J_{\pm}^{\theta}. \end{eqnarray}

\smallskip

{\bf Example. 3}. When r=2, i.e.
\begin{eqnarray}
\theta=\pmatrix{
0&0&0 \cr
1&0&0 \cr
0&1&0 \cr}\;\;\;\;\;\;\;\hbox{and}\;\;\;\;\;\;\;J_{\pm}^{\theta}=\pmatrix{
J_{\pm}^{(0)}\;&\;\;0\;\;&0 \cr
\delta\;J_{\pm}^{(1)}&J_{\pm}^{(0)}&0 \cr
\delta^{2} J_{\pm}^{(2)}&\delta\;J_{\pm}^{(1)}&J_{\pm}^{(0)} \cr}, 
\nonumber
\end{eqnarray}
we obtain
\begin{eqnarray}
&&[H,\;J_{\pm}^{\theta}]=\pm 2 J_{\pm}^{\theta}, \nonumber\\ 
&&[J_{+}^{\theta},\;J_{-}^{\theta}]=H+\theta^{2} {\delta^{3} \over 3!}
(H^{3}-H). 
\end{eqnarray}
The coproduct is given by
\begin{eqnarray}
&&\bigtriangleup (H) = H \otimes \hbox{{\bf 1}} + 
\hbox{{\bf 1}} \otimes H, \\
&&\bigtriangleup(J_{\pm}^{\theta})=J_{\pm}^{\theta}\;\bar\otimes\; 
(\hbox{{\bf 1}}+{\theta \delta\over 2}H+{(\theta\delta)^{2} 
\over 8}H^{2})+(\hbox{{\bf 1}}-{\theta\delta\over 2}H+
{(\theta\delta)^{2} \over 8}H^{2})\;\bar\otimes\;J_{\pm}^{\theta}. 
\nonumber 
\end{eqnarray}
Such a structure is discussed in $[7]$ in connection 
with the Higgs algebra, characterized by 
\begin{eqnarray}
&&[H,\;J_{\pm}]=\pm 2 J_{\pm}, \nonumber\\ 
&&[J_{+},\;J_{-}]=H+c\; 
H^{3}, 
\end{eqnarray}
$c$ being an arbitrary constant. This algebra is of special 
interest as it is the one of dynamical symmetries for the harmonic 
ascillator and the Kepler problem in a two-dimensional curved space $[8]$.

\smallskip

{\bf Example .4.} The $r\rightarrow \infty$ case ($\theta$ is 
equivalent to a {\bf real} variable) is characterized by
\begin{eqnarray}
&&J_{\pm}^{\theta}=\sum_{m=0}^{\infty}J_{\pm}^{(m)}\;\delta^{m}
\;\theta^{m} \nonumber\\
&&\phantom{J_{\pm}^{\theta}}=\sum_{m=0}^{\infty}J_{\pm}^{(m)}\;\zeta^{m} 
\nonumber\\
&&J_{\pm}^{\theta}:={\tilde J}_{\pm} 
\end{eqnarray}
and
\begin{eqnarray}
&&[H,\; {\tilde J}_{\pm}]=\pm 2 {\tilde J}_{\pm}, \nonumber\\
&&[{\tilde J}_{+},\;{\tilde J}_{-}]={e^{\zeta H}-e^{-\zeta H} \over
e^{\zeta}-e^{-\zeta}}, 
\end{eqnarray}
where $\zeta=\theta\;\delta$. We thus recover the Drinfeld-Jimbo structure 
${\cal U}_{\displaystyle{e^{\zeta}}}(sl(2))$ as a particular case of 
${\cal U}^{\infty}_{\displaystyle{q}}(sl(2))$.
 
The same embedding is also present at the level of the Hopf structure
with
\begin{eqnarray}
&&\bigtriangleup (H) = H \otimes \hbox{{\bf 1}} + \hbox{{\bf 1}} \otimes H,
\nonumber \\
&&\bigtriangleup ({\tilde J}_{\pm}) = {\tilde J}_{\pm}\otimes e^{\zeta H/2} + 
e^{-\zeta H/2}\otimes {\tilde J}_{\pm},\nonumber \\
&& \varepsilon ({\tilde J}_{\pm})=\varepsilon (H)=0, \nonumber\\
&& S(H)=-H,\;\;S({\tilde J}_{\pm})=-e^{{\pm} \zeta}\;{\tilde J}_{\pm}. 
\end{eqnarray}

\sect{The ${\cal U}_{q_1,q_2}^{r_{1},r_{2}}(sl(2))$ algebra}

Let us now introduce, for example, two real paragrassmannian variables 
$\theta_{1}$ and $\theta_{2}$ respectively of order $r_{1}$ and $r_{2}$, i.e.
\begin{eqnarray}
&&\theta_{1}^{r_{1}+1}=0,\;\;\;\;\;\;\;\;\;\theta_{2}^{r_{2}+1}=0,
\nonumber \\
&&\theta_1 \theta_2 +\theta_2 \theta_1 = 0. 
\end{eqnarray}

Using the Campbell-Baker-Hausdorff expansion
\begin{eqnarray}
(\exp A)(\exp B)=\exp C,
 \end{eqnarray}
where
\begin{eqnarray}
&&C=A +B + {1\over 2} \sum_{m=1}^{\infty}
{1 \over (m+1)!}(ad\;A)^{m}(B)+{1\over 2} \sum_{m=1}^{\infty}
{(-1)^{m} \over (m+1)!}(ad\;B)^{m}(A),\nonumber \\
&&(ad\;A)^{m}(B)=\lbrack A,\lbrack A,\cdots , \lbrack A,\;B\rbrack \rbrack 
\cdots \rbrack, \nonumber\\
&&(ad\;B)^{m}(A)=\lbrack B,\lbrack B,\cdots , \lbrack B,\;A\rbrack \rbrack 
\cdots \rbrack, 
\end{eqnarray}
we propose to define
\begin{eqnarray}
J_{\pm}^{(\theta_{1}, \theta_{2})}=\sum_{m=0}^{\infty} \theta^{m}J_{\pm}^{(m)}
\end{eqnarray}
where
\begin{eqnarray}
&&\theta = \theta_1 \delta_1 + \theta_2 \delta_2+ {1 \over 2} 
\sum_{m=1}^{r_{1}}
{\delta_2\delta_1^{m} 2^{m} \over (m+1)!} \theta_{1}^{m}\theta_2 + 
{1 \over 2} \sum_{m=1}^{r_{2}}
{\delta_1 \delta_2^{m} 2^{m} \over (m+1)!} \theta_{1}\theta_{2}^{m}\nonumber \\
&&\exp (\theta_1 \delta_1 ) \exp ( \theta_2 \delta_2)= \exp \theta.  
\end{eqnarray}
Using (4.4), we deduce that 
\begin{eqnarray}
&&[J_{+}^{(\theta_{1},\theta_{2})},\;J_{-}^{(\theta_{1},\theta_{2})}]=
{e(H\delta_1;\;\theta_{1})e(H\delta_2;\;\theta_{2})-
e(-H\delta_2;\;\theta_{2}) e(-H\delta_1;\;\theta_{1}) \over
e(\delta_1;\;\theta_{1})e(\delta_2;\;\theta_{2})-
e(-\delta_2;\;\theta_{2}) e(-\delta_1;\;\theta_{1})},\nonumber\\
&&[H,\;J_{\pm}^{(\theta_{1},\theta_{2})}]=
\pm 2 J_{\pm}^{(\theta_{1},\theta_{2})},
\label{eq:2-theta-algebra} 
\end{eqnarray}
The algebra $\{ J_{\pm}^{(\theta_{1},\theta_{2})},H\}$ described by 
the commutation relations (\ref{eq:2-theta-algebra}) is just the    
quantization of $sl(2)$ with two paragrassmannian variables and is 
denoted by ${\cal U}_{\delta_1,\delta_2}^{\theta_1 , \theta_2 }(sl(2))$. The 
${\cal U}_{\delta_1,\delta_2}^{\theta_1 , \theta_2 }(sl(2))$ algebra is 
equipped with the following Hopf structure
\begin{eqnarray}
&&\bigtriangleup(H)=H\otimes\hbox{{\bf 1}}+\hbox{{\bf 1}} 
\otimes H, \nonumber\\
&&\bigtriangleup (J_{\pm}^{(\theta_1,\theta_2 ) }) = 
J_{\pm}^{(\theta_1,\theta_2 )}\;\bar \otimes \; 
e\;({H\delta_1\over 2};\;\theta_1)\;e\;({H\delta_2\over 2};\;\theta_2) +  
e\;(-{H\delta_2\over 2};\;\theta_2)\;e\;(-{H\delta_1\over 2};\;\theta_1) 
\;\bar \otimes \; J_{\pm}^{(\theta_1,\theta_2 )},\nonumber \\
&& \varepsilon (J_{\pm}^{\theta})=\varepsilon (H)=0,\qquad\qquad
\varepsilon(\hbox{{\bf 1}})=\hbox{{\bf 1}}, \nonumber\\
&& S(H)=-H,\qquad\qquad S(\hbox{{\bf 1}})=\hbox{\bf 1}, \\
&&S(J_{\pm}^{(\theta_1,\theta_2 )})= 
-e\;({H\delta_1\over 2};\;\theta_1)\;e\;({H\delta_2\over 2};\;\theta_2)
J_{\pm}^{(\theta_1,\theta_2 )}
e\;(-{H\delta_2\over 2};\;\theta_2)\;e\;(-{H\delta_1\over 2};\;\theta_1), 
\nonumber\\  
&&\bigtriangleup(e\;({H\delta_1\over 2};\;\theta_1)\;e\;({H\delta_2
\over 2};\;\theta_2))= 
e\;({H\delta_1\over 2};\;\theta_1)\;e\;({H\delta_2\over 2};\;\theta_2)
\;\bar \otimes \; e\;({H\delta_1\over 2};\;\theta_1)\;
e\;({H\delta_2\over 2};\;\theta_2).\nonumber 
\end{eqnarray}

\sect{Connection with some Nonlinear Algebras} 

Let us take in ${\cal U}_{q}^{r}(sl(2))$ the following 
restriction
\begin{eqnarray}
q=1\qquad\qquad\hbox{i.e.}\qquad\qquad
q=e^{2\pi i n}, \qquad\qquad 
\end{eqnarray}
where $n$ characterizes the Riemann branch. We have
\begin{eqnarray}
e\;(2\pi i  n \Omega ;\;\theta)
=\cos(2\pi n \Omega ;\;\theta)+i \sin(2\pi n \Omega ;\;\theta), 
\end{eqnarray}
with
\begin{eqnarray}
&&\cos(x ;\;\theta)=\sum_{k=0}^{[r/2]} (-1)^{k} 
{x^{2k} \theta^{2 k} \over (2 k )!}, \nonumber\\
&&\sin(x ;\;\theta)=\sum_{k=0}^{(r-1)/2-{1 \over 2}(1+(-1)^{r})} (-1)^{k} 
{x^{2k+1} \theta^{2 k+1} \over (2 k+1 )!}. 
\end{eqnarray}
Thus, the commutation relations are written as
\begin{eqnarray}
&&[H,\;J_{\pm}^{\theta}]=\pm 2 J_{\pm}^{\theta}, \\
&&[J_{+}^{\theta},\;J_{-}^{\theta}]={\sin(2 \pi n H;\; \theta) \over
\sin(2 \pi n ;\; \theta)}.\nonumber
\end{eqnarray}

When $n\rightarrow \infty$, we deduce
\begin{eqnarray}
\lim_{n\rightarrow \infty}{\sin(2 \pi n H;\; \theta) \over
\sin(2 \pi n ;\; \theta)} = H^{r-{1 \over 2}(1+(-1)^{r})} 
\end{eqnarray}
and
\begin{eqnarray}
&&[H,\;J_{\pm}^{\theta}]=\pm 2 J_{\pm}^{\theta}, \nonumber\\
&&[J_{+}^{\theta},\;J_{-}^{\theta}]=H^{r-{1 \over 2}(1+(-1)^{r})},  
\end{eqnarray}
the deformation is being a nonlinear one. 

\smallskip

Now, if we take in ${\cal U}_{q_1,q_2}^{r_{1},r_{2}}(sl(2))$ 
$r_{2}\longrightarrow \infty$ and $\delta_1=2\pi i n$ 
($n\longrightarrow \infty$), we deduce the following nonlinear algebra
\begin{eqnarray}
&& [J_{+},J_{-}]={H^{r_{1}}q^{H}-(-1)^{r_{1}}q^{-H}H^{r_{1}}
\over q-(-1)^{r_{1}}q^{-1}}, \nonumber\\
&&[H, J_{\pm}]=\pm J_{\pm}. 
\end{eqnarray}

\sect{ Conclusion}

We have proposed new deformed structures ${\cal U}_{q}^{r}(sl(2))$ 
and ${\cal U}_{q_1,q_2}^{r_{1},r_{2}}(sl(2))$ obtained by 
paragrassmannian deformation. When the order of the paragrassmannian
variable goes to infinity, we recover the 
Drinfeld-Jimbo scheme of deformation.

\smallskip

It has also to be noticed that our proposal points out two different 
Hopf structures for the same deformed algebra. In particular, $sl(2)$ 
can be associated with a cocommutative coproduct ($r=0$) or a 
non-cocommutative one ($r=1$). Then it is possible to get a new 
${\cal R}$-matrix given by
\begin{eqnarray}
&&{\cal R}_{\theta}=1\otimes 1 +\delta\theta(J_{-}\otimes J_{+}-
J_{+}\otimes J_{-}) \nonumber \\
&& \phantom{ {\cal R}_{\theta} }=U_{\theta}U_{-\theta}^{+},
\end{eqnarray}
where 
\begin{eqnarray}
U_{\theta}=1\otimes 1 +{1\over 2}\delta\theta(J_{-}\otimes J_{+}-
J_{+}\otimes J_{-}) 
\end{eqnarray}
by requiring 
\begin{eqnarray}
U_{\theta}\bigtriangleup_{r=0} (a)
= \bigtriangleup_{r=1} (a) U_{\theta},
\end{eqnarray}
for any a belonging to $sl(2)$. It has also to be noticed that this matrix 
${\cal R}_{\theta}$ satisfies the Yang-Baxter equation. Thus it is the first 
solution, to our knowledge, depending on a paragrassmannian variable.

\smallskip

We would like to mention that the $r=2$-case is a particulary
interesting one as already mentionned. It is the first case where the 
deformation is present at the level of the algebra and these deformations 
are nonlinear ones in the sense of Ro${\tilde c}$ek. We have thus defined 
ad-hoc coproducts, counits and antipodes for such deformations being of 
physical interest. 

\smallskip

Finally, the restriction of the values of the parameters of the 
deformation gives somes nonlinear algebras as particular cases.

\vskip 1cm

{\bf Acknowledgments}  

We thank Daniel Arnaudon and Jean Lascoux for important discussions and 
encouragement.

Two of us (J.B. and N.D.) also thank the two other authors (B.A. and A.C.) for 
their warm hospitality during their stay in Palaiseau where part of this 
work has been elaborated under a TOURNESOL financial support 
(which is also acknowledged by all of us).


\vfill
\vfill
\vfill

\newpage

\begin{thebibliography}{99}


\bibitem {}
V.G.Drinfeld, Quantum Groups, Proc. Int. Congress of Mathematicians, Berkeley,
California, Vol. {\bf 1}, Academic Press, New York (1986) 798.

M.Jimbo, Lett. Math. Phys. {\bf 10} (1985) 63.

\bibitem {}
V. Chari and A. Pressley, A Guide to Quantums Groups, Cambridge (cambridge UP)
(1994). 

\bibitem {}
L.C. Biedenharn, J. Phys. A {\bf 22} (1989), L873.

A.J. Macfarlane, J. Phys. A {\bf 22} (1989), 4581.

P. Roche and D. Arnaudon, Lett. Math. Phys. {\bf 17} (1989), 295.

\bibitem {}
M. Ro${\tilde c}$ek, Phys. Lett. B {\bf 255} (1991), 554.

\bibitem {}
C. Quesne, Phys. Lett. A {\bf 293} (1994), 245.

\bibitem {}
Y. Ohnuki and S. Kamefuchi, J. Math. Phys. {\bf 24} (1980), 609.

\bibitem {} B. Abdesselam, J. Beckers, A. Chakrabarti and N. Debergh,
On Nonlinear Angular Momentum Theories, Their Representations 
and Associated Hopf Structures, To be published in J. Phys. A: Math. Gen. 
(1996). 


\bibitem {}
P.W. Higgs, J. Phys. A {\bf 12} (1979), 309.



\end{thebibliography}
\end{document}